%% file: LandscapeFormation.tex
\documentclass[aps,prx,reprint,superscriptaddress,longbibliography]{revtex4-1}

\usepackage{amssymb}
\usepackage{amsmath}
\usepackage{bm}
\usepackage{graphicx}
\usepackage{xcolor}
\usepackage{hyperref}

\newtheorem{lemma}{Lemma}
\newtheorem{remark}{Remark}

\DeclareMathOperator{\M}{M}

\begin{document}

\title{Construction of quasi-potentials for stochastic dynamical systems: an optimization approach}

\author{R. D. Brackston}
\email[]{r.brackston13@imperial.ac.uk}
\affiliation{Department of Life Sciences, Imperial College London, SW7 2AZ, United Kingdom}

\author{A. Wynn}
\affiliation{Department of Aeronautics, Imperial College London, SW7 2AZ, United Kingdom}

\author{M. P. H. Stumpf}
\affiliation{Department of Life Sciences, Imperial College London, SW7 2AZ, United Kingdom}

\date{\today}

\begin{abstract}
The construction of effective and informative landscapes for stochastic dynamical systems has proven a long-standing and complex problem. In many situations, the dynamics may be described by a Langevin equation while constructing a landscape comes down to obtaining the quasi-potential, a scalar function that quantifies the likelihood of reaching each point in the state-space. 
In this work we provide a novel method for constructing such landscapes by extending a tool from control theory: the Sum-of-Squares method for generating Lyapunov functions. Applicable to any system described by polynomials, this method provides an analytical polynomial expression for the potential landscape, in which the coefficients of the polynomial are obtained via a convex optimization problem. The resulting landscapes are based upon a decomposition of the deterministic dynamics of the original system, formed in terms of the gradient of the potential and a remaining ``curl'' component. By satisfying the condition that the inner product of the gradient of the potential and the remaining dynamics is everywhere negative, our derived landscapes provide both upper and lower bounds on the true quasi-potential; these bounds becoming tight if the decomposition is orthogonal. The method is demonstrated to correctly compute the quasi-potential for high-dimensional linear systems and also for a number of nonlinear examples.
\end{abstract}

\maketitle

\section{Introduction}    \label{sec:intro}

Multi-dimensional nonlinear dynamical systems can exhibit complex behavior including limit cycles, strange attractors and multiple fixed points. When also driven by stochastic perturbations, such systems often explore a range of conditions and may exhibit behavior such as random switching between attractors and stochastic resonance \citep{Gammaitoni1998}. The mathematics of such systems has its roots in the description of Brownian motion, which motivated the study of stochastic differential equations, generally written in the form,
\begin{equation}    \label{eq:sde}
	\text{d}\bm{x} = f(\bm{x})\,\text{d}t + g(\bm{x})\,\text{d}W_t.
\end{equation}
Here the deterministic component of the dynamics (the drift) is described by $f(\bm{x})$, while the stochastic component (diffusion) is given by $g(\bm{x})\,\text{d}W_t$, where $\text{d}W_t$ describes the increment of a Wiener process \cite{Gardiner:2009aa}.

A popular description of such systems is in terms of an energy landscape, in which the dynamics are described as a ball moving in a potential basin. In some contexts, the landscape may simply give an intuitive description of the dynamics \cite{Rigas2015,Brackston2016}, while in others it may represent a true energy function \cite{Wolynes1995}. In the context of developmental biology, Waddington's epigenetic landscapes provide a popular analogy for stem cell development \citep{Waddington1957,Moris2016}, even though it has proven impossible to define a true energy function for general developmental processes.

If correctly formed, landscapes may offer a quantitative analysis and interrogation tool beyond a merely phenomenological description, even in cases where the free energy of the system may not be defined. In particular, a so-called quasi-potential landscape may be defined in terms of the Freidlin-Wentzell action functional \cite{Freidlin2012}. This landscape provides quantitatively a measure of the relative stability of different states and the most probable paths between them. The concept of landscapes defined by the action has received particular recent attention \cite{Wang2010b}, as has the development of methods to calculate the so-called Minimum Action Path (MAP) \cite{Heymann2008,Liu2008,Delacruz2018}. Other recent methods have sought to evaluate quasi-potential landscapes with a focus on the action functional. For linear systems an analytic expression for the quasi-potential was obtained by \cite{Chen2005}, while numerical methods applicable to nonlinear systems have been given in \cite{Cameron2012,Dahiya2017,Cameron2017}. While such methods may offer a numerical solution over a discretized space, in this work we present a method that generates an analytical solution, and furthermore is applicable to both linear and nonlinear systems.

Given the motivation behind landscape descriptions, it is unsurprising that considerable effort has been applied to find other methods for landscape evaluation based either upon experimental data or mathematical models \cite{Wang2015,Zhou2016}. In biological applications, the most popular and readily applied method is that based upon the steady-state probability distribution \cite{Wang2008}. This approach is based upon the Fokker-Planck equation, which justifies computing a potential as $U(\bm{x}) \propto -\ln(P_S(\bm{x}))$.
In practice the steady-state distribution is either obtained by solving the Fokker-Planck equation directly \cite{Lv2014,Ge2015,Li2016}, or found from extensive simulations of the corresponding SDE \cite{Wang2010a,Li2011,Guo2017}. While this method is easy to implement, it may be impractical for higher-dimensional systems and there is generally no guarantee that the derived landscape relates directly to the fixed points of the deterministic system.

Irrespective of how a landscape is obtained, its existence implies a decomposition of the deterministic vector field into two components, one of which is given by the gradient of the landscape while the other accounts for the remainder of the dynamics in $f(\bm{x})$:
\begin{equation} \label{eq:Decomposition}
	f(\bm{x}) = -\nabla U(\bm{x}) + f_U(\bm{x}).
\end{equation}
The construction of the landscape may therefore also be performed by considering the properties of this decomposition directly. As discussed by \cite{Freidlin2012} and advocated by \cite{Zhou2012}, one particularly interesting case is that in which the gradient and remainder are everywhere orthogonal, implying an independence between gradient-based and rotational dynamics. While this property may be satisfied in some cases, it is important to note that it remains unclear if this is always possible. Motivated by this, we develop here a method to form landscapes satisfying the slightly more general \emph{sub-orthogonal} condition that $\nabla U(\bm{x}) \cdot f_U(\bm{x}) \leq 0$, of which the orthogonal decomposition is a special case.

As we shall discuss below, the sub-orthogonal decomposition has the desirable property that under certain conditions it is directly related to the quasi-potential. One of the other desirable properties of the landscape, and one that is satisfied by the sub-orthogonal decomposition, is that it should be a Lyapunov function \cite{Jost:2005} for the deterministic dynamics. In this context, Lyapunov functions are often used to prove the asymptotic stability of the fixed points $f(\bm{x})=0$. However, they also have the useful property that they provide the basins of attraction for the deterministic dynamics. There is therefore a strong equivalence between quasi-potential landscapes and Lyapunov functions, suggesting that pre-existing methods for generating such functions may also be used to obtain the quasi-potential. This is the approach that we take here.

In this paper, we develop a method that performs a sub-orthogonal decomposition of the deterministic dynamics $f(\bm{x})$, by utilizing a popular method for the construction of Lyapunov functions: the Sum-of-Squares (SOS) method \cite{Parillo2000}. We will firstly provide some background on the properties of the sub-orthogonal decomposition in section~\ref{sec:background} before providing details of this computational method in section~\ref{sec:methods}. We will then examine its application to a series of examples. The first of these examples will be linear systems in section~\ref{sec:linear} for which reference analytical solutions exist, followed by application to nonlinear systems with various properties in section~\ref{sec:nonlinear}. Conclusions will finally be given in section~\ref{sec:conclusions}.

\section{Properties of the sub-orthogonal decomposition}    \label{sec:background}

In this work we consider the case in which the deterministic part of \eqref{eq:sde} is decomposed as in \eqref{eq:Decomposition}, and this decomposition has the property that,
\begin{equation} \label{eq:Suborthogonality}
	\nabla U(\bm{x}) \cdot f_U(\bm{x}) \leq 0.
\end{equation}
When this relation becomes an equality, the decomposition can be said to be orthogonal, while otherwise we will refer to the decomposition as \emph{sub-orthogonal}. As we shall now discuss, such a decomposition of the vector field has a useful relation to the quasi-potential and furthermore provides a Lyapunov function for the deterministic dynamics.

\subsection{The quasi-potential}

Let us now restrict our analysis to systems in which the diffusion tensor $g(x)$ is a uniform diagonal matrix, describing purely additive noise,
\begin{equation}    \label{eq:sde2}
	\text{d}\bm{x} = f(\bm{x})\,\text{d}t + \sigma I_n\,\text{d}W_t.
\end{equation}
Here, $\sigma$ is a small constant and $I_n$ is the $n \times n$ identity matrix. Given such a system, the action associated with a path $\bm{x}(t) = \varphi$ is given by the following integral,
\begin{equation}	\label{eq:S}
	S(\varphi) = \frac{1}{\sigma^2} \int_0^T \left|\dot{\varphi}-f(\varphi)\right|^2 \,\text{d}t.
\end{equation}

Equation~\eqref{eq:S} is useful because in the limit,  $\sigma \to 0$, the probability of the stochastic system following a certain path is directly related to the action of that path according to,
\begin{equation}
	P(\varphi) \propto e^{-S(\varphi)}.
\end{equation}

Given the properties of the action, one can define a quasi-potential $Q_a$ with respect to a fixed point $a$ as,
\begin{equation}
	Q_a(\bm{x}) = \inf_{\varphi,T} \left[S(\varphi)|\varphi(0)=a, \varphi(T)=\bm{x}\right].
\end{equation}
This quasi-potential therefore gives the action associated with the minimum action path to each point within the basin of attraction of the fixed point, and in turn may be related to the probability distribution over the state space.

Given the sub-orthogonality condition of \eqref{eq:Suborthogonality}, the quasi-potential may be evaluated as,
\begin{align*}
	Q_a(\bm{x}) & = \frac{1}{\sigma^2} \inf_{\varphi,T} \left[ \int_0^T \left| \dot{\bm{x}}+\nabla U-f_U \right|^2 \,\text{d}t\right] \\
    & = \frac{1}{\sigma^2} \inf_{\varphi,T} \left[ \int_0^T \left| \dot{\bm{x}}-\nabla U-f_U \right|^2 + \right. \nonumber\\
    & \qquad\qquad\qquad\qquad  4\dot{\bm{x}}\cdot\nabla U - 4f_U\cdot\nabla U \,\text{d}t \Bigg] \\
    & = \frac{4}{\sigma^2} (U(\bm{x})-U(a)) + \nonumber \\
    &  \frac{1}{\sigma^2} \inf_{\varphi,T} \left[ \int_0^T \underbrace{\left| \dot{\bm{x}}-\nabla U-f_U \right|^2}_{\text{(i)}} - \underbrace{4f_U\cdot\nabla U}_{\text{(ii)}} \,\text{d}t \right]
\end{align*}
In the orthogonal case, term (ii) is equal to zero, while the infimum of term (i) will approach zero for a path arbitrarily close to that satisfying the ODE $\dot{\bm{x}} = \nabla U + f_U$ \footnote{Note that because $a$ is a fixed point, such a path will take infinite time. However given any small perturbation away from $a$, the time will be finite.}. For the sub-orthogonal case the net contributions of both terms will be positive. The potential $U$ therefore provides a lower bound for the quasi-potential according to,
\begin{equation}
	Q_a(\bm{x}) \geq \frac{4}{\sigma^2} (U(\bm{x})-U(a)).
\end{equation}
In the case of the orthogonal decomposition this inequality becomes an equality.

\subsection{Lyapunov functions}

Lyapunov functions serve an important purpose in the field of non-linear dynamical systems and feedback control, in proving asymptotic stability of a fixed point. While the proof of asymptotic stability requires specific conditions to be satisfied, here we use a slightly looser definition that allows the Lyapunov function to have wider interpretation. For a system defined by a set of ODEs $\dot{\bm{x}}=f(\bm{x})$ a Lyapunov function $U(\bm{x}): \mathbb{R}^n \to \mathbb{R}$, is one that satisfies the following constraints:
\begin{subequations}    \label{eq:lyapunov}
\begin{align}
	U(\bm{x}) \geq 0 & , \label{eq:lyapunovA} \\
    \nabla U(\bm{x}) \cdot f(\bm{x}) \leq 0 & . \label{eq:lyapunovB}
\end{align}
\end{subequations}
These constraints impose the requirements that $U$ is everywhere positive and that $U$ decreases along trajectories of $f(\bm{x})$. The key implication of these constraints is that the Lyapunov function correctly captures the basins of attraction for all the stable fixed points: a region of the state space $\mathcal{A}$ around a fixed point $a$ such that if $\bm{x}(0) \in \mathcal{A}$, $\bm{x}(t) \to a$ as $t \to \infty$. Such a function therefore provides an accurate qualitative description of a landscape underlying the dynamics.

Given the sub-orthogonal decomposition property of \eqref{eq:Suborthogonality}, we may rewrite constraint \eqref{eq:lyapunovB} as,
\begin{align}
	\nabla U \cdot f & = \nabla U \cdot (-\nabla U + f_U) \nonumber \\
    				 & = -|\nabla U|^2 + \nabla U \cdot f_U \leq 0.
\end{align}
Sub-orthogonality of the decomposition is therefore a sufficient condition for the potential $U$ to be a Lyapunov function for the deterministic dynamics described by $f$. We therefore use this condition in our optimization approach discussed below.

\section{The Optimization Method}    \label{sec:methods}

In this work we develop a method for landscape formation based on the sum-of-squares method for forming Lyapunov functions. We will therefore first give a brief overview of the standard algorithm through which SOS is used to find Lyapunov functions, before detailing the key changes required in our algorithm to achieve a sub-orthogonal decomposition. The algorithm itself has been implemented in the programming language Julia, built upon the JuMP package for mathematical optimization \citep{Dunning2017} and in \textsc{Matlab} using the SOStools package \cite{sostools}.

\subsection{The standard SOS method}

Finding a Lyapunov function for nonlinear dynamical systems is generally a difficult problem\cite{Jost:2005,Silk:2011dc}, and is known to be NP-hard in the case where $U$ is a polynomial \cite{Murty1987}: i.e. given a possible function $U$, even checking that it satisfies the requirements may be very computationally expensive. The basis of the SOS method is that polynomial functions formed as the sum of the squares of lower order polynomials are guaranteed to be positive-definite, and furthermore can be found via efficient optimization methods. While the requirement that $U(\bm{x})$ is SOS is stricter than that of positive-definiteness, it makes the problem computationally tractable.

The standard use of SOS in forming a Lyapunov function involves solving the following:
\begin{subequations}    \label{eq:SOSproblem}
    \begin{align}
	& \text{find a feasible} & & [c_j]_{j = 0}^m \nonumber \\
    & \text{subject to} & & U = c_0 + \sum_{j=1}^m c_j p_j \label{eq:SOSbasis} \\
    & & & U \geq \epsilon \sum_{i=1}^n x_i^2, \label{eq:SOSbndCon} \\
    & & & \nabla U \cdot f \leq 0. \label{eq:SOSdecCon}
    \end{align}
\end{subequations}
Here the $p_j$ are a suitable set of $m$ monomial terms (e.g. $x_1x_2^2$) and the $c_j$ are real coefficients. The parameter $\epsilon$ is a positive constant, typically of order one.
The task therefore involves a pure feasibility problem, and comes down to finding any $U$ composed of the monomials specified in \eqref{eq:SOSbasis}, subject to the inequality constraints \eqref{eq:SOSbndCon}, \eqref{eq:SOSdecCon} \footnote{In practice, these inequalities are imposed by specifying that $U(\bm{x}) - \epsilon \sum_i x_i^2 \in \Sigma$ and $-\nabla U(\bm{x}) \cdot f(\bm{x}) \in \Sigma$, where $\Sigma$ denotes the set of polynomials that are Sum of Squares.}. There are therefore infinitely many possible solutions, if $f(\bm{x})$ describes a stable system. As we will discuss below, we modify each of the steps of problem \eqref{eq:SOSproblem} to produce a convex optimization problem that yields a unique result.

\subsection{Towards orthogonality}    \label{subsec:orthogonality}

While the SOS method is able to generate Lyapunov functions for general $n$-dimensional dynamical systems, such functions are not unique and will generally not satisfy the sub-orthogonality requirement of \eqref{eq:Suborthogonality}. As noted by several previous authors \cite{Zhou2012,Cameron2012}, an orthogonality requirement results in a Hamilton-Jacobi equation $\nabla U\cdot f + ||\nabla U||^2 = 0$ which is a nonlinear constraint on $U(\bm{x})$. Such constraints cannot be directly implemented into the SOS program. In order to implement such a constraint in a linear manner, consider the matrix
\begin{equation}
	M_U := \begin{bmatrix}
    		  -\nabla U \cdot f & \nabla U^\top \\
    		   \nabla U & I_n \end{bmatrix} \in \mathbb{R}^{(n+1)\times(n+1)}.
\end{equation}
A Schur Complement argument \cite{Zhang2005} implies that $M_U \succeq 0$ is positive definite for any $\bm{x} \in \mathbb{R}^n$ if and only if
\begin{equation}
	-\nabla	U \cdot f - \nabla U^\top I_n \nabla U \geq 0,
\end{equation}
which is equivalent to $\nabla U \cdot f_U \leq 0$, where $f_U$ is as defined in \eqref{eq:Decomposition}. This in turn implies that $\nabla U \cdot f \leq 0$, which satisfies the key requirement for a Lyapunov function and constraint \eqref{eq:SOSdecCon}.
With the matrix inequality constraint in place, there will still be an infinite number of possible $U$, most of which will not come close to orthogonality. In order to find the unique $U$ that minimizes $|\nabla U \cdot f_U|$, we find it necessary to make $U$ as steep as possible, in this way maximizing the contribution of $-\nabla U$ to $f$. This is achieved by maximizing an appropriately chosen lower bound,
\begin{equation}    \label{eq:bound}
	B(\bm{x})=\sum_i\epsilon_i b_i(\bm{x}),
\end{equation}
where each $b_i$ is a monomial in $\bm{x}$ and the $\epsilon_i$ are real coefficients. The choice of these monomials is important and the method for choosing them is described below in section~\ref{subsec:monomial}. The complete optimization procedure is therefore as follows:

\begin{subequations}    \label{eq:Orthproblem}
    \begin{align}
	& \underset{[c_j], [\epsilon_i]}{\text{maximize}} & & \sum_i \epsilon_i \label{eq:objective} \\
    & \text{subject to} & & U = c_0 + \sum_{j=1}^m c_j p_j \label{eq:basis} \\
    & & & U \geq \sum_i \epsilon_i b_i, \label{eq:bndCon} \\
    & & &	M_U \succeq 0. \label{eq:MvCon}
    \end{align}
\end{subequations}

For the case that $f(\bm{x})$ is linear, this method can be proven to yield the correct result, as detailed in appendix~\ref{app:proof1}.

\subsection{Choosing the monomial basis and lower bound}    \label{subsec:monomial}

The first key choice that must be made is that of the monomial basis from which $U$ is constructed, i.e. the $p_j$ in \eqref{eq:basis}. In the typical application of SOS to form Lyapunov functions, a candidate function is formed from a full collection of the monomials in $\bm{x}$ up to a given (even) degree $d$ \cite{Papachristodoulou2002}. For a given $n$-dimensional system, this results in $m = \binom{n+d}{d}$ monomials forming the polynomial function $U$ \cite{Majumdar2014}. For example, for a 4th order polynomial and 4-dimensional state space this results in 40 monomials. Since the computational cost of the optimization increases with $m$, it is advantageous to reduce the basis if possible, by exploiting particular properties of the function we are trying to obtain.

\subsubsection{A minimal basis}

We can initially provide an upper bound on the degree $d$ using the sub-orthogonality condition \eqref{eq:Suborthogonality}, which may be rewritten as,
\begin{equation}
	\nabla U \cdot f + |\nabla U|^2 \leq 0.
\end{equation}
Since $|\nabla U|^2 \geq 0$, the above inequality can only hold if the degree of $\nabla V \cdot f$ is at least as big as the degree of $|\nabla U|^2$. If $d$ is the degree of $U$ and $e$ the degree of $f$, then $d + e - 1 \geq 2d - 2$, meaning that,
\begin{equation}    \label{eq:maxdegree}
	d \leq e + 1.
\end{equation}

While \eqref{eq:maxdegree} provides an upper bound on the degree of $U$, we may further refine the basis \eqref{eq:basis}, motivated by the properties in the pure gradient case, in which $f_i = -\partial U/\partial x_i$. The minimal basis for such a $U$ is therefore provided by
\begin{equation}    \label{eq:minbasis}
	[p_j]_{j=1}^m = \M \left[ \sum_i f_i x_i \right],
\end{equation}
where the operator $\M$ extracts a vector of monomials from a polynomial.

As an example, consider the vector field defined by,
\begin{equation}
	f(\bm{x}) = \begin{bmatrix} x_1 - x_1^3 \\
												 -x_2^3	 \\
												 -x_3^3 \end{bmatrix}.
\end{equation}
Here the highest order in $f$ is three so \eqref{eq:maxdegree} implies that we will likely require a $d=4$ order polynomial. Given that the system has $n=3$ dimensions, the standard method leads to 35 monomials. By contrast, a minimal basis would only include five terms, $[x_1^4,x_2^4,x_3^4,x_1^2,1]$.

\subsubsection{The lower bound}

Given a monomial basis for $U$, it is next necessary to choose the monomials $b_i$ in the lower bound \eqref{eq:bound}. This choice is best illustrated via the following simple system described by,
\begin{align}	\label{eq:eg1D}
	f(x) & = x - x^3 + \beta \nonumber \\
		 & = -\frac{\text{d}}{\text{d}x}\left[\frac{1}{4}x^4 - \frac{1}{2}x^2 - \beta x + C \right].
\end{align}
Given that \eqref{eq:eg1D} is one-dimensional and therefore may be written as a pure gradient system, the minimal basis argument of \eqref{eq:minbasis} correctly ascertains that $U(x)$ will be composed only of monomials $[x^4,x^2,x,1]$. A sensible choice for the lower bound will be $b(x)=x^q$, where $q$ is an even integer. The inequality constraint~\eqref{eq:bndCon} therefore becomes:
\begin{equation}
	c_3x^4 + c_2x^2 + c_1x + c_0 \geq \epsilon x^q,
\end{equation}
where the choice of $q$ is of critical importance. In particular it must satisfy two particular properties:

\paragraph{Sufficient pressure}
Suppose that we take $q=2$, now for any $c_3,c_0>0$ the constraint may be satisfied. Yet $c_3$ may remain arbitrarily small and the problem will not have converged. The lower bound is therefore unable to exert sufficient ``upward pressure'' on $U$.

\paragraph{Feasibility}
Alternatively consider $q=6$. In this case it is impossible to find coefficients such that $U>\epsilon x^6$, since as $x \to \infty$, $c_3x^4 < \epsilon x^6$ for any $c_3,\epsilon>0$. The program will therefore be infeasible.

\begin{figure}
    \centering
    \includegraphics[width=0.45\textwidth]{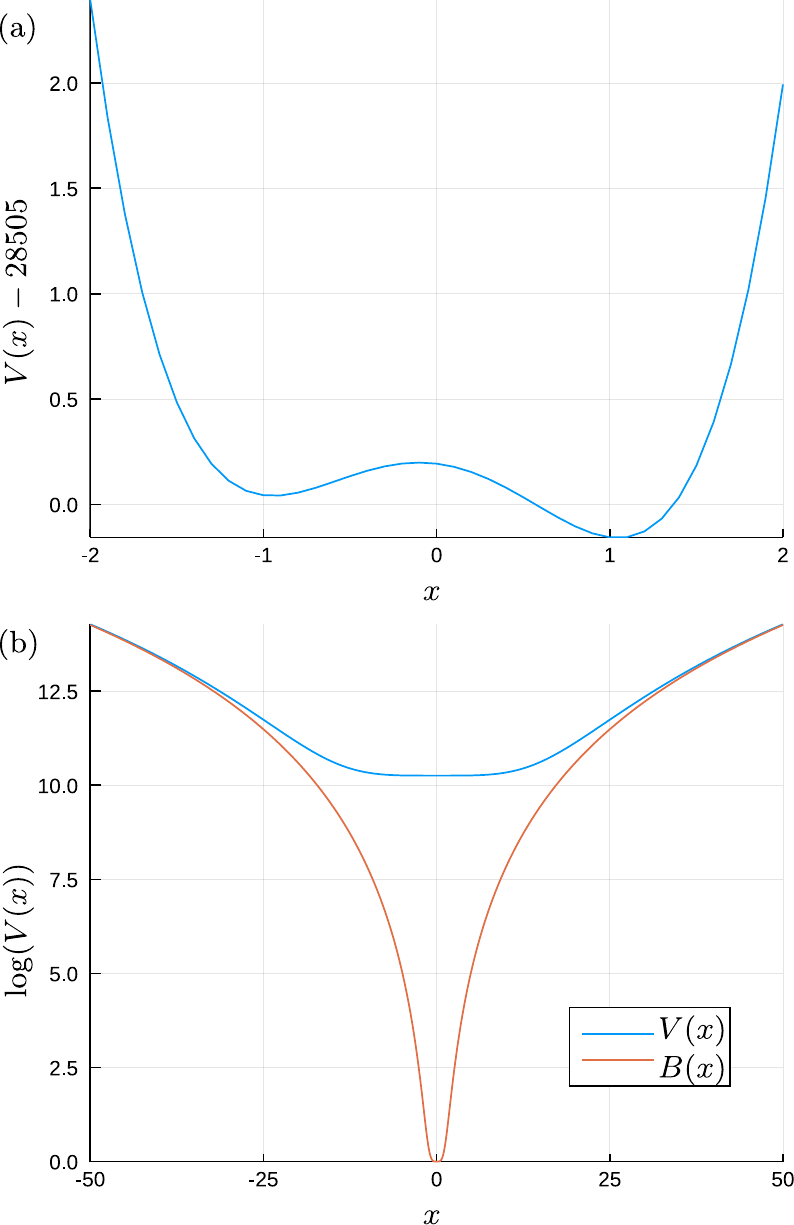}
    \caption{A schematic of the lower bounding method. (a) The potential function  in the region of the fixed points and (b) a comparison between the potential and lower bound.}
   \label{fig:BoundSchematic}
\end{figure}

The correct choice for $q$ in this case is $q=4$, i.e. the highest degree in the basis for $U$. An example of this lower bound in practice is displayed in figure~\ref{fig:BoundSchematic}. The same principles of sufficient upward pressure and feasibility also apply to higher dimensional systems, such that the lower bound is chosen to consist of the highest degree single and mixed monomials in the basis for $U$.

\subsubsection{Extending the basis}    \label{subsubsec:extbasis}

A correctly chosen ``minimal'' basis will naturally work for the case that $f(\bm{x})$ is a pure gradient system, and may sometimes work in other cases \footnote{For a linear system in which the matrix $A$ is normal, a minimal basis will also be sufficient since sparsity of the matrix is maintained through \eqref{eq:linearU}}. However more generally, our constraints may require that $\nabla U$ contains monomials beyond \eqref{eq:minbasis}, as will be illustrated in later examples. Nonetheless, in addition to considerations of computational cost, the sensitivity of the method to the choice of lower bound means that we cannot simply choose all monomials up to degree $d$, as is done in the standard SOS method.

Consider again the example given in \eqref{eq:eg1D}. Based on the above discussion it is clear that we could not have included an additional higher-order monomial in the basis for $U$. For example had we included a $x^6$ term, we would logically have also chosen $q=6$ for the lower bound. Since the actual potential need not include a $x^6$ term, this would make the program infeasible.
The principles of sufficient pressure and feasibility in the choice of the lower bound, therefore guide which monomial terms are permitted in the basis for $U$. The procedure for determining the basis is therefore as follows:
\begin{enumerate}
\item The highest order monomial in any individual $x_i$ should be that obtained from the minimal basis.
\item Other mixed monomials are required, provided that:
	\begin{enumerate}
	\item the total order is not more than the highest total order in the minimal basis,
	\item the individual order in each $x_i$ is not more than the highest individual order monomial in the minimal basis.
	\end{enumerate}
\end{enumerate}

\subsection{Iterative improvement}    \label{subsec:iterative}

The process outlined above performs well for many systems, but may require further improvement to push the obtained landscape closer to orthogonality. Given an initial Lyapunov function $U_1$ such as that found by implementing \eqref{eq:Orthproblem}, one can attempt to iteratively improve the solution, generating at each step a second Lyapunov function $U_2$ composed of the same monomial basis. Each iteration involves solving the following optimization problem:
\begin{subequations} \label{eq:iterProblem}
    \begin{align}
	& \underset{[c_j], \alpha, \epsilon}{\text{minimize}} & & \alpha \label{eq:objective2} \\
    & \text{subject to} & & U_2= c_0 + \sum_{j=1}^m c_j p_j \label{eq:basis2} \\
    & & & U_2 \geq \epsilon \sum_i x_i^2, \label{eq:posDef} \\
    & & &	M_{U_2} \succeq 0, \label{eq:Mv2} \\
    & & & \nabla U_2 \cdot(f+2\nabla U_1) \geq \nonumber \\
    & & & \quad\,\, \alpha(f\cdot \nabla U_1) + (1+\alpha)||\nabla U_1||^2, \label{eq:alpha} \\
    & & & \alpha,\epsilon > 0.
    \end{align}
\end{subequations}
The key to this optimization is that by using the initial guess, $U_1$, the optimization is constructed to be linear in the new improved Lyapunov function, $U_2$.
Constraints~\eqref{eq:posDef} and \eqref{eq:Mv2}, simply mirror those in optimization~\eqref{eq:Orthproblem}, namely positive definiteness of the Lyapunov function and the sub-orthogonality condition. Constraint~\eqref{eq:alpha} is chosen such that for $\alpha=1$ equality is guaranteed if $U_2=U_1$, while for $\alpha < 1$, $\nabla U_2 \cdot f_{U_2} \leq \nabla U_1 \cdot f_{U_1}$. This result is proven in Appendix~\ref{app:proof}, however the key point is that the optimization has a guaranteed feasible solution with $\alpha=1$ and $U_2=U_1$, while if a smaller $\alpha$ is obtained then the new Lyapunov function is closer to orthogonality.

The complete method for obtaining the potential landscape therefore consists of first applying procedure~\eqref{eq:Orthproblem}, followed by a number of iterations of procedure~\eqref{eq:iterProblem}. We generally find that only one or two iterations of this method are required for satisfactory convergence. The number of iterations may therefore either be pre-specified or determined based upon a convergence criterion.

\section{The linear case}    \label{sec:linear}

A special case of the orthogonal decomposition may be considered when the deterministic system is linear \cite{Cameron2012}. While this may seem a very limiting test case, the linear situation can exemplify some of the key scenarios and limitations that also apply in the nonlinear case, while remaining tractable to analysis.

\subsection{Properties of linear systems}

Linear dynamics may be written in terms of a matrix multiplication as,
\begin{equation}
	\dot{\bm{x}} = f(\bm{x}) = A\bm{x},
\end{equation}
where the matrix $A \in \mathbb{R}^{n \times n}$. The construction of the quasi-potential then comes down to decomposing the matrix $A$ into two parts as
\begin{equation}
	A = A_g + A_c,
\end{equation}
where the gradient matrix $A_g$ is a symmetric matrix that gives the potential according to
\begin{equation}
	U(\bm{x}) = -\frac{1}{2} \bm{x}^\top A_g \bm{x}.
\end{equation}
Because $A_g$ is symmetric, it has purely real eigenvalues, while if the decomposition is orthogonal, the remainder matrix $A_c$ has purely imaginary eigenvalues and is such that the product $A_g A_c$ is an antisymmetric matrix.

In such linear cases it has been demonstrated that an orthogonal decomposition does always exist \cite{Kwon2005}, and furthermore may be provided by an analytical expression \cite{Chen2005} as,
\begin{equation}
	A_g = \frac{1}{2} \left( \int_0^\infty e^{At}e^{A^*t}\,\text{d}t  \right)^{-1}.
\end{equation}
If the matrix $A$ is normal ($AA^* = A^*A$), then this expression simplifies to,
\begin{equation}    \label{eq:linearU}
	A_g = \frac{1}{2} \left( A+A^*\right) \bm{x}.
\end{equation}
Since these expressions provide solutions against which our method can be tested, linear systems give a perfect test case for our optimization-based approach which does not rely on any prior knowledge of the properties of the system.

\subsection{A three-dimensional example}

We consider now the test-case of a three-dimensional linear system defined by,
\begin{align}
	f(\bm{x}) & = A\bm{x} \nonumber \\
					& = \begin{bmatrix} -5.0 & 0.0 & 0.2 \\
										0.0 & -1.5 & 3.0 \\
										0.5 & -5.0 & -1.0 \end{bmatrix} \bm{x}.
\end{align}
Implementation of methods \eqref{eq:Orthproblem} and \eqref{eq:iterProblem} provides the following decomposition for $A$:
\begin{multline}
	A = \underbrace{\begin{bmatrix} -5.01 & 0.14 & 0.18 \\
						   0.14 & -1.55 & -0.02 \\
						   0.18 & -0.02 & -0.94 \end{bmatrix}}_{A_g} + \\
        \underbrace{\begin{bmatrix} 0.01 & -0.14 & 0.02 \\
						   -0.14 & 0.05 & 3.02 \\
						   0.32 & -4.98 & -0.06 \end{bmatrix}}_{A_c},
\end{multline}
giving the quasi-potential as,
\begin{multline}
	U(\bm{x}) = 2.50x_1^2 + 0.78x_2^2 + 0.47x_3^2 - \\
    			0.14x_1x_2 - 0.18x_1x_3 + 0.02x_2x_3.
\end{multline}
As can readily be observed, the matrix $A_g$ is symmetric and can easily be verified to have purely real (negative) eigenvalues. The matrix $A_c$ is not itself either symmetric or antisymmetric but can be verified to have purely imaginary eigenvalues. The dynamics associated with $A_c$ are therefore those of pure oscillation, and evolve on level sets of the potential $U$.

It is noteworthy that for the system defined by $A$, there is no direct connection between $x_1$ and $x_2$. One may therefore naively expect such a direct connection to also be absent in the potential. Nevertheless, in order to provide an orthogonal decomposition, the potential matrix $A_g$ does include such terms that are then negated in $A_c$. This confirms that a truly minimal basis is indeed insufficient for forming the quasi-potential, as discussed in section~\ref{subsubsec:extbasis}.

\subsection{Scaling with system size}

\begin{figure}
    \centering
    \includegraphics[width=0.4\textwidth]{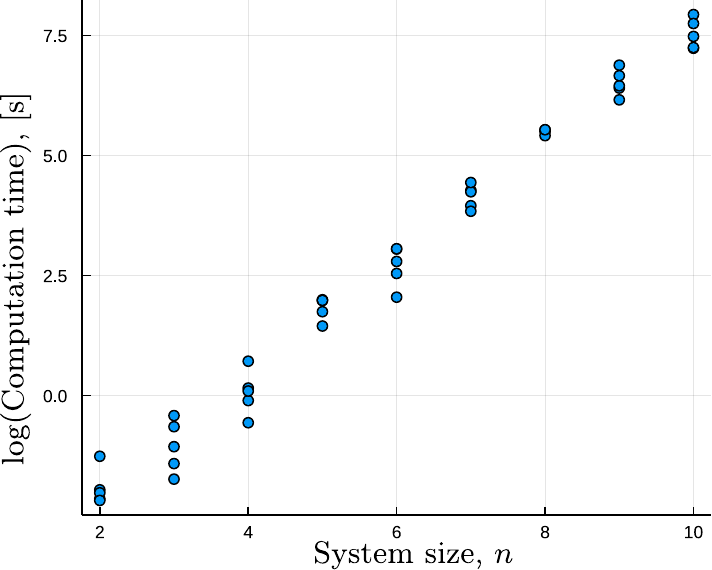}
    \caption{Scaling of the computational cost with system size. Tests were performed on a basic desktop PC with 16~Gb ram and an Intel Core i3 processor.}
   \label{fig:Scaling}
\end{figure}

While the method may be verified to work for this three-dimensional example, it is also useful to assess the accuracy and efficiency of the algorithm for higher dimensional systems. We do this by randomly generating real, negative-definite matrices of varying size $n = 2:10$. In all such cases the method is able to perform a decomposition that satisfies the orthogonality constraint. Results for the computational cost are shown in figure~\ref{fig:Scaling}, displayed on a logarithmic scale. It is evident that the cost increases considerably with system size, most likely in a factorial fashion. For small systems the algorithm converges in less than a second, while for larger systems it may take up to an hour.

It is worth noting here that some of the additional computational cost associated with larger systems arises from a greater number of times that the iteration procedure \eqref{eq:iterProblem} must be applied, given a desired level of orthogonality. For low-dimensional problems we find that sometimes no iterations are required, while for systems of dimension $n=10$, up to five iterations may be necessary to achieve the same quality of result.

\section{The nonlinear case}    \label{sec:nonlinear}

Having studied the ability of our method to obtain the decomposition in linear cases, we now move on to fully non-linear examples for which analytical solutions do not always exist.

\subsection{A nonlinear multistable system}

Our first example is the quartic system from \cite{Zhou2012} with four attractors and a known orthogonal decomposition. The dynamics are defined by,
\begin{equation}
	f(\bm{x}) = \begin{bmatrix} -1 + 9x_1 - 2x_1^3 + 9x_2 - 2x_2^3 \\
												   1 - 11x_1 + 2x_1^3 + 11x_2 - 2x_2^3  \end{bmatrix},
\end{equation}
for which the true potential is given by,
\begin{equation}
	U(\bm{x}) = 0.5(x_1^4  + x_2^4) - 5(x_1^2 + x_2^2) + x_1x_2 + x_1.
\end{equation}
For such a problem our method finds the potential to within 5 significant figures within a few seconds. Because the output of the algorithm is a symbolic expression for the potential, the actual value of $U$ may then be evaluated at any required points in the space with minimal further effort. This is in contrast of course, to methods that solve the associated PDE over a discretized grid, since evaluation outside of the solved area requires significant further computation.

\subsection{The Maier-Stein model}

We next apply the method to the widely studied Maier-Stein model of \cite{Maier1996}, defined by,
\begin{equation}
	f(\bm{x}) = \begin{bmatrix} x_1 - x_1^3 - \gamma x_1x_2^2 \\
										-\mu(1+x_1^2)x_2  \end{bmatrix}	.
\end{equation}
This model has the property that when $\gamma=\mu$, the dynamics can be expressed in terms of a pure potential \footnote{This may be verified by checking that $\frac{\partial f_1}{\partial x_2}$ = $\frac{\partial f_2}{\partial x_1}$.} $f(\bm{x})=-\nabla U$, where,
\begin{equation}
	U(\bm{x}) = \underbrace{+0.25}_{c_1}x_1^4 \underbrace{+0.5\mu}_{c_2}x_1^2x_2^2  \underbrace{-0.5}_{c_3}x_1^2 \underbrace{+0.5\mu}_{c_4}x_2^2.
\end{equation}
For cases in which $\gamma \neq \mu$, the system also has a non-gradient component.

We apply our method for three different parameter combinations, obtaining in each case a polynomial expression for the potential with the same non-zero terms. The coefficients for these terms are displayed in table~\ref{tab:MaierStein}. For the non-gradient cases, the decompositions are, as desired, sub-orthogonal.

\begin{table}
\caption{\label{tab:MaierStein}Coefficients obtained for the Maier-Stein model for three different parameter combinations.}
\begin{ruledtabular}
\begin{tabular}{cccccc}
 $\mu$ & $\gamma$ & $c_1$ & $c_2$ & $c_3$ & $c_4$\\ \hline
 1.0 & 1.0 & 0.2499999 & 0.4999999 & - 0.5000001 & 0.4999998 \\
 1.0 & 10.0 & 0.1600785 & 1.000557 & - 0.3199390 & 0.0003017 \\
 2.5 & 1.0 & 0.2495550 & 0.5109047 & - 0.4991115 & 1.249649 \\
\end{tabular}
\end{ruledtabular}
\end{table}

\subsection{Bounds on the quasi-potential}

A useful property of the (sub-)orthogonal decomposition is that it may be used to obtain predicted MAPs. As discussed in section~\ref{sec:background}, MAPs are those paths that minimize the action functional and provide a definition of the true quasi-potential. In general these paths must be found via a costly optimization that must be performed for every point in the state-space. For a system that satisfies the orthogonal decomposition however, the MAP from a fixed point $\bm{x}_o$ to another point $\bm{x}_e$ within the basin of attraction, is one that follows,
\begin{equation}    \label{eq:predMAP}
	\dot{\bm{x}} = \nabla U + f_U.
\end{equation}
Such a path may therefore readily be obtained by a simulation starting at $\bm{x}_e$ and following $\dot{\bm{x}} = -\nabla U - f_U$.

\begin{figure*}
    \centering
    \includegraphics[width=0.9\textwidth]{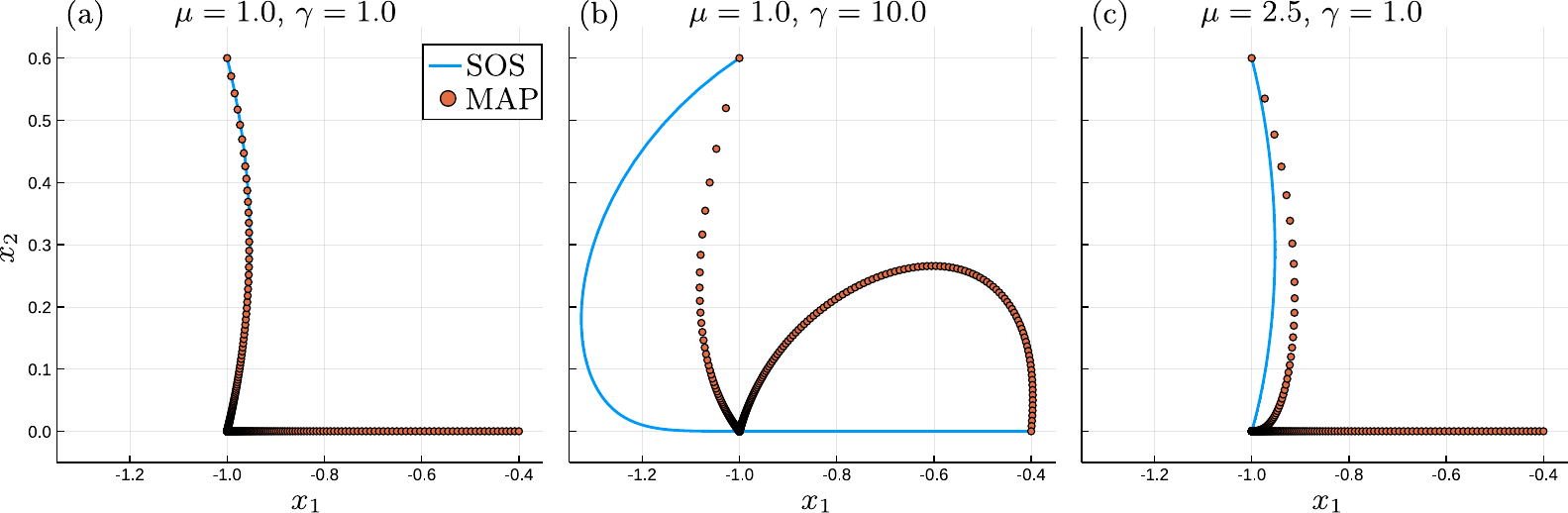}
    \caption{A comparison of the SOS estimated and true MAPs. Paths start at the fixed point $\bm{x}_o=(-1.0,0.0)$ and end at the points $\bm{x}_e=(-0.4,0.0)$ and $\bm{x}_e=(-1.0,0.6)$. Points along the minimum action paths are equally spaced in time.}
   \label{fig:Paths}
\end{figure*}

While we expect paths described by \eqref{eq:predMAP} to match exactly in the orthogonal case, for sub-orthogonal cases we hope that there may still be close agreement, depending on the degree of orthogonality. We therefore compare the predicted and exact paths, evaluating the true MAPs via our own implementation of the optimization approach detailed in \cite{Delacruz2018}. A comparison of paths from the fixed point $\bm{x}_o=(-1.0,0.0)$ is given in figure~\ref{fig:Paths}. In the gradient case, the paths match exactly, as expected. For the first non-gradient cases, there is clearly some disagreement, although the sense of curvature of the paths is generally the same for the predicted and true MAPs.

While the agreement of the paths may seem quantitatively poor, these paths may still be used to provide an estimate of the true quasi-potential. Given that the quasi-potential is defined in terms of the infimum of the action over all possible paths, the action of any path that we choose can provide an upper bound to the true quasi-potential, evaluated via the geometric minimum actio method \cite{Heymann2008}. While most choices of possible path will give an action much greater than that of the MAP, the SOS predicted path may be expected to give a tighter bound, owing to the qualitative similarity with the true MAP. In addition to this upper bound on the true quasi-potential, the potential from the decomposition itself provides a lower bound, as discussed in section~\ref{sec:background}.

\begin{figure*}
    \centering
    \includegraphics[width=0.9\textwidth]{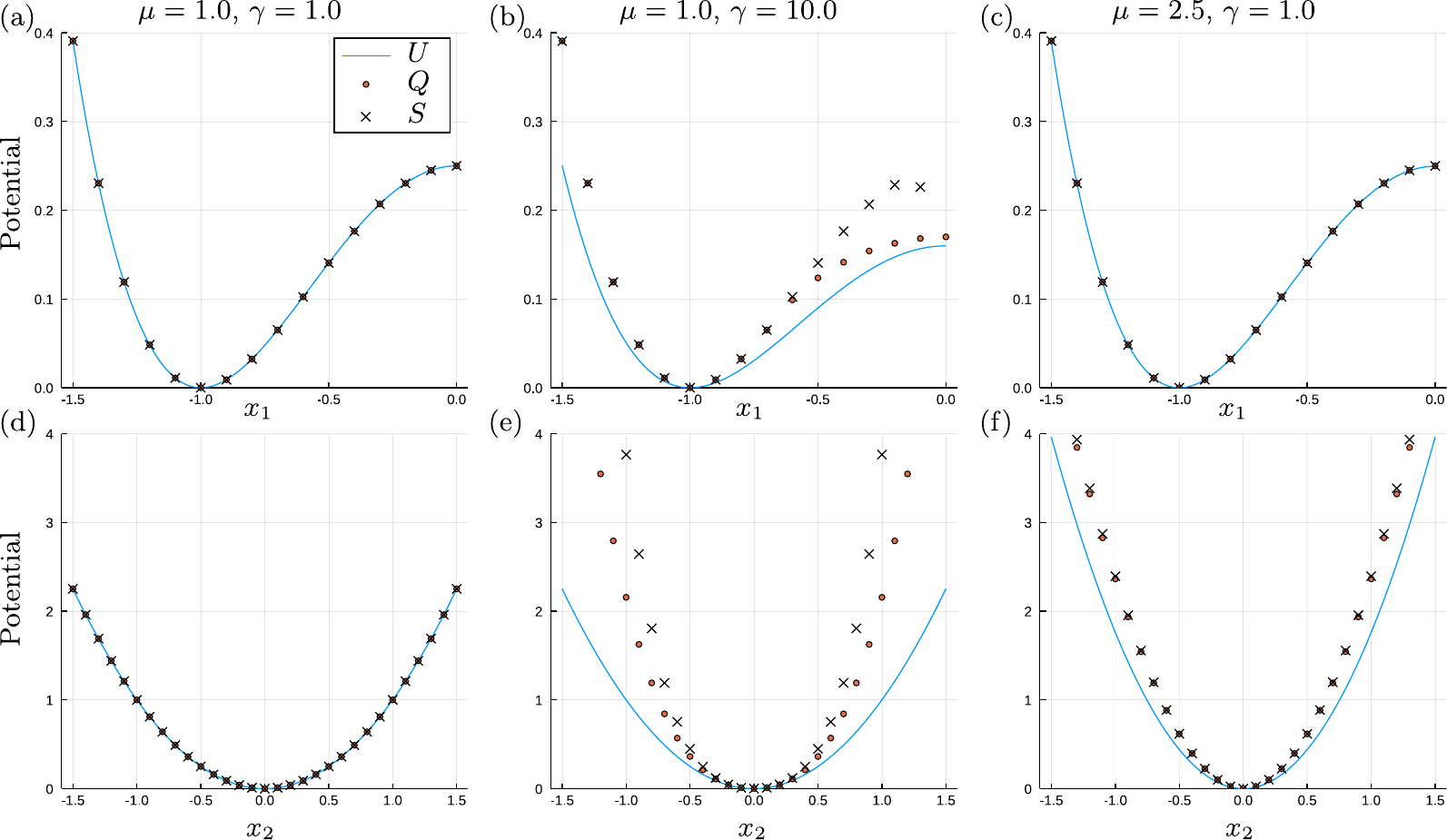}
    \caption{A comparison of the estimates of the quasi-potential. $U$ is the potential formed from the sub-orthogonal decomposition, $Q$ is the true quasi-potential and $S$ is the action of the predicted minimum action path.}
   \label{fig:Bounds}
\end{figure*}

A quantitative comparison between the true quasi-potential $Q$, the potential from the decomposition $U$ and that from the SOS predicted path $S$ is given in figure~\ref{fig:Bounds}. In the gradient case there is almost exact agreement between all three methods, demonstrating that the obtained decomposition provides both the true quasi-potential and predicts the correct paths. In the first non-gradient case, the upper bound remains remarkably tight, despite the large differences between the true and predicted paths. For the second non-gradient case the bounds are even tighter, with an exact match along the $x$-axis. In all cases, each of $U$ and $S$ can be seen to indeed be lower and upper bounds respectively.

\section{Conclusions}    \label{sec:conclusions}

Evaluating the quasi-potential for linear and nonlinear SDEs is a challenging problem, yet one with significant interest and motivation. In this work we have provided a novel method for the calculation of the quasi-potential based upon the Sum of Squares method for constructing Lyapunov functions. Our method is applicable to systems for which the governing equations are polynomial and involves solving an optimization over the coefficients of a polynomial potential function.

The construction of an informative landscape is motivated by three key requirements. Firstly, we would like the potential to correctly capture the basins of attraction for the deterministic system. Such a requirement is equivalent to the potential being a Lyapunov function, and is therefore naturally achieved by our method.
Secondly, it is desirable to have an estimate of the most probable transition trajectories between basins of attraction; the so-called minimum action paths. For cases permitting an orthogonal decomposition of the dynamics, the paths may readily be obtained from the two vector-field components. For cases in which the decomposition is sub-orthogonal, the obtained paths may provide a more-or-less accurate approximation, suitable for use as an initial guess in an optimization routine.
Finally, we may wish for the potential function to be a quasi-potential for the system, accurately describing the transition probabilities for situations of vanishingly small noise. Again, for systems permitting an orthogonal decomposition of the dynamics our method calculates exactly such a quasi-potential, while for sub-orthogonal cases the potential may be used to provide both a lower and upper bound.

The first key limitation of the method is in its applicability to only polynomial systems. However, such systems are commonplace in e.g. mass action models of chemical kinetics and linear models for dynamical systems. A closely related limitation is in the ability of the method to express the potential itself in terms of a polynomial. It is plausible that in some cases even if the governing equations are polynomial, a potential satisfying the normal decomposition must be expressed in some expanded basis beyond monomial terms. Regardless, our polynomial sub-orthogonal potential still provides useful bounds.

A second key limitation is that the obtained quasi-potential is only valid for systems with additive noise, in which the noise tensor is equal to the identity matrix. Yet this is a common approximation, especially when the magnitude of the noise tends to zero, and furthermore can always be achieved for linear systems via a coordinate transform. If multiplicative noise is present but may expressed as a polynomial, it is possible that in some cases this could be incorporated into the algorithm, however this is beyond the scope of this study.

While we have provided a method that generates an analytical expression for the quasi-potential for a particular subset of SDE systems, it remains unclear if such a feat can be achieved in more general cases. Ultimately, we hope that the method provides another useful tool for the interpretation and analysis of stochastic dynamical systems, and may pave the way for more general methods to obtain the quasi-potential in a symbolic manner.

\begin{acknowledgments}
We are grateful to Dr David Schnoerr for critical reading of the manuscript. This work was supported by BBSRC grant BB/N003608/1.
\end{acknowledgments}

\section*{Code availability}

The algorithms discussed in this work have been implemented in both Julia and \textsc{Matlab}, freely accessible at \url{https://github.com/rdbrackston/SDEtools} and \url{https://github.com/rdbrackston/normalSoS} respectively.

\appendix

\section{Optimality in the linear case} \label{app:proof1}

In section~\ref{subsec:orthogonality}, we presented a method by which we maximize a lower bound for the quasi-potential $U$, attempting to make the potential as steep as possible. Here we justify this method by proving that in the linear case, the matrix defining a normal decomposition is also that which provides the largest potential, defined in a suitable way.

\begin{lemma}
Consider a full rank real matrix $A \in \mathbb{R}^{n \times n}$ whose eigenvalues have negative real part, and $U(\bm{x}) = \frac{1}{2}\bm{x}^\top P \bm{x}$ for some $P=P^\top \in \mathbb{R}^{n \times n}$. \\
Then the decomposition,
\begin{equation*}
	A\bm{x} = -\nabla U(\bm{x}) + f_U(\bm{x}), \qquad \bm{x} \in \mathbb{R}^{n},
\end{equation*}
holds with $f_U(\bm{x}) = (A+P)\bm{x}$. \\
Suppose further that $P$ solves the optimization problem:

\begin{subequations} \label{eq:Optimisation}
    \begin{align}
	& \underset{P}{\text{maximize}} & & \text{tr} (P) \label{eq:objective} \\
    & \text{subject to} & & P = P^\top, \label{eq:constraint1} \\
    & & & P(A+P) \preceq 0. \label{eq:constraint2}
    \end{align}
\end{subequations}
Then $PA+A^\top P + 2P^2 = 0$. Hence, $\nabla U \cdot f_U = 0$.
\end{lemma}

\noindent \textit{Proof}. First note that \eqref{eq:Optimisation} is feasible since $A$ is stable. Furthermore, for any $P$ satisfying $P(A+P) \leq 0$, we have,
\begin{equation*}    \label{eq:Cdef}
	C := 2P^2 + PA + A^\top P \leq 0.
\end{equation*}
Hence, the Ricatti equation,
\begin{equation*}    \label{eq:ricatti}
	 2P^2 - P(-A) -(-A)^\top P - C = 0,
\end{equation*}
has a (trivial) solution $P=P^\top \in \mathbb{R}^{n \times n}$. Now since $C'=0 \geq C$, it follows from \cite[][Theorem 2.3]{Gohberg1986}, that the Riccati equation,
\begin{equation*}    \label{eq:ricatti2}
\begin{split}
	 & 2X^2 - X(-A) - (-A)^\top X = \\
     & 2X^2 - X(-A) - (-A)^\top X - C' = 0,
\end{split}
\end{equation*}
has a maximal solution $X=P_+$, which satisfies $P_+ \succeq P$.
So, we have proved that there exists a solution $P_+$ to,
\begin{equation*}
	2P_+^2 + P_+A + A^\top P_+ = 0,
\end{equation*}
which has the property that $P_+ \succeq P$ for any feasible variable of \eqref{eq:Optimisation}. Let $\hat{P}$ be the optimal decision variable of \eqref{eq:Optimisation}. Then, since $P_+$ is also feasible for \eqref{eq:Optimisation}, 
\begin{equation}    \label{eq:trace}
	\text{tr}(\hat{P}) \geq \text{tr}(P_+).
\end{equation}
However,
\begin{align} \label{eq:pdef}
	P_+ \succeq \hat{P} & \Rightarrow (P_+)_{ii} \geq (\hat{P})_{ii} \nonumber \\
    					& \Rightarrow (P_+ - \hat{P})_{ii} \geq 0. 
\end{align}
By \eqref{eq:trace},
\begin{align*}
   & & \text{tr}(P_+ - \hat{P}) \leq & 0 \nonumber \\
   & \Rightarrow & \sum_{i=1}^n \underbrace{(P_+ - \hat{P})_{ii}}_{\text{all $\geq 0$ by \eqref{eq:pdef}}} \leq & 0 \\
   & \Rightarrow & (P_+)_{ii} = (\hat{P})_{ii}, & \qquad i=1,...,n.
\end{align*}
Hence, $P_+ \succeq \hat{P}$ and both have the same diagonal entries. Therefore $P_+ = \hat{P}$. \\
Consequently, the solution $\hat{P}$ to \eqref{eq:Optimisation} is the maximal solution to $\hat{P}A+A^\top \hat{P} + 2\hat{P}^2 = 0$
\hfill $\square$

\begin{remark}
In practice, for linear systems our full method \eqref{eq:Orthproblem} actually implements,
\begin{subequations}  \label{eq:RealOptimisation}
    \begin{align}
	& \underset{P,c_0,[\epsilon_i]}{\text{maximize}} & & \sum_{i=1}^n \epsilon_i \label{eq:Realobjective} \\
    & \text{subject to} & & \frac{1}{2}\bm{x}^\top P \bm{x} + c_0 \geq \sum_{i=1}^n \epsilon_i x_i^2, \label{eq:Realconstraint1} \\
    & & & P(A+P) \preceq 0. \label{eq:Realconstraint2}
    \end{align}
\end{subequations}
While \eqref{eq:RealOptimisation} this is not entirely equivalent to \eqref{eq:Optimisation}, \eqref{eq:Realconstraint1} $\Rightarrow \sum_i \epsilon_i \leq \text{tr}(P)$. The implemented method therefore provides a very close approximation that also generalizes to the nonlinear case.

\end{remark}

\section{Iterative improvement} \label{app:proof}

In section~\ref{subsec:iterative} we displayed the iterative optimization approach through which an improved Lyapunov function could be obtained, given a suitable initial guess. The optimization procedure may be justified by proving the following lemma.

\begin{lemma}
Let $U : \mathbb{R}^n \to \mathbb{R}$ be fixed. Suppose that there exists $\alpha \geq 0$ and $V : \mathbb{R}^n \to \mathbb{R}$ such that
\begin{subequations}    \label{eq:IterCond}
\begin{align}
	&V \geq 0, \label{eq:posDefA} \\
    &M_{V} \succeq 0, \label{eq:MvA} \\
    &\nabla V \cdot (f+2\nabla U) \geq \alpha(f\cdot \nabla U) + (1+\alpha)||\nabla U||^2 \label{eq:alphA}
\end{align}
\end{subequations}
Then
\begin{equation*}
	\alpha(\nabla U \cdot f_{U}) + ||\nabla V - \nabla U||^2 \leq \nabla V \cdot f_{V} \leq 0. 
\end{equation*}
\end{lemma}

\noindent \textit{Proof}. First lower bound $\nabla V \cdot f_{V}$:
\begin{equation*}
\begin{split}
	\nabla V \cdot f_{V} & = \nabla V \cdot (f+\nabla V) \\
    & = \nabla V \cdot (f+2\nabla U) + ||\nabla	V||^2 - 2\nabla V\cdot \nabla U \\
    (\text{by \eqref{eq:alphA}}) & \geq \alpha(f\cdot \nabla U + ||\nabla U||^2) + ||\nabla U||^2 + \\
     & \qquad\qquad\qquad\qquad\quad\,\,\,\, ||\nabla V||^2 - 2\nabla V \cdot \nabla U \\
    & = \alpha(\nabla U\cdot f_U) + ||\nabla U - \nabla V||^2.
\end{split}
\end{equation*}
That $\nabla V \cdot f_V \leq 0$ follows from \eqref{eq:MvA}. \hfill $\square$
 
\begin{remark}
(i) If $U$ satisfies the standard Lyapunov conditions $U \geq 0$, $M_U \succeq 0$ then the choice ($\alpha=1$, $V=U$) also satisfies conditions \eqref{eq:IterCond}. Hence, there is always a feasible point which satisfies these inequalities.

(ii) If ($\alpha$, $V$) satisfiers the constraints with $0 < \alpha < 1$, then
\begin{equation*}
	\alpha(\nabla U \cdot f_{U}) + ||\nabla V - \nabla U||^2 \leq \nabla V \cdot f_{V} \leq 0
\end{equation*}
indicates that an improved lower-bound has been found for $V \cdot f_V$. The improvement over $U \cdot f_U$ is determined by the difference $||\nabla V - \nabla U||^2$ and the scaling factor $\alpha$.

(iii) If $\alpha=0$ and $V$ satisfies the constraints, then \eqref{eq:alphA} reads
\begin{equation*}
	||\nabla U - \nabla V||^2 \leq \nabla V \cdot f_V \leq 0
\end{equation*}
which forces $V=U$ and indicates that $U \cdot f_U=0$, i.e. the original Lyapunov function was a good choice.
\end{remark}

\input{LandscapeFormation.bbl}

\end{document}

%% file: LandscapeFormation.bbl
%